\documentclass[10pt]{article}
\newtheorem{theorem}{Theorem}
\newtheorem{lemma}{Lemma}
\newtheorem{remark}{Remark}
\newtheorem{corollary}{Corollary}
\def\Frac#1#2{\frac{\displaystyle{#1}}{\displaystyle{#2}}}
\usepackage{amssymb}
\usepackage{epsfig}

\begin{document}

 \title{The Schwarzian-Newton method for solving nonlinear equations, with applications}

\author{
Javier Segura \thanks{The author acknowledges financial support from Ministerio de Econom\'{\i}a
y Competitividad (project MTM2012-34787)}\\
        Departamento de Matem\'aticas, Estad\'{\i}stica y 
        Computaci\'on,\\
        Univ. de Cantabria, 39005 Santander, Spain.\\
        e-mail: javier.segura@unican.es
}

\date{\ }

\maketitle
\begin{abstract}
The Schwarzian-Newton method can be defined as the minimal method for solving nonlinear
equations $f(x)=0$ which is exact for any function $f$ with constant Schwarzian derivative; exactness means that 
the method gives the exact root in one iteration for 
any starting value in a neighborhood of the root. 
This is a fourth order method which has Halley's method 
as limit when the Schwarzian derivative tends to zero. 
We obtain conditions for the convergence of the SNM in an interval and show how this method can be 
applied for a reliable and fast solution of some problems, like the inversion of cumulative
distribution functions (gamma and beta distributions) and the inversion of elliptic integrals. 
\end{abstract}

2000MSC: Primary 65H05; Secondary 33B20, 33E05.

\section{Introduction}

The problem of inverting functions is one of the central problems in numerical analysis, with uncountable
applications. There exists a vast literature on numerical methods for solving nonlinear equations
and many different methods with different properties are available; particularly, there
is a considerable interest in building high order methods with the highest possible efficiency index \cite{Traub:1964:IMF}.

The study of non-local convergence properties is, however, not 
so common. But without a knowledge of these properties, the use of fast numerical methods 
will require accurate initial estimations in order to guarantee convergence. There
is, in addition, the general rule that the higher the convergence rate, the more
unpredictable the method will be when accurate initial approximations are not available.
This is probably one of the reasons why high order methods are not so common in applications; 
bisection, secant and Newton's methods
are the most popular followed closely by Halley's method (order $3$) and other third order variants. 

The ideal situation would be that of a high order method with good and manageable non-local convergence properties; however,
there is little (if any) information on the non-local behaviour of methods of order larger than three, particularly concerning
verifiable conditions for convergence in an interval.
In this article, we propose the Schwarzian-Newton (SNM) which can defined as the minimal method for solving nonlinear
equations $f(x)=0$ which is exact for any function $f$ with constant Schwarzian derivative, meaning that the method gives 
the exact root in one iteration for 
any starting value in a neighborhood of the root. 
This is a fourth order method which has Halley's method (HM) as limit when the Schwarzian derivative tends to zero. 
As happens for the HM \cite{Melman:1997:GAC,Salehov:1952:OTC}, the SNM has a simple geometrical interpretation in terms of osculating curves.

The SNM will be proven to converge (monotonically) in less iterations that Halley's method (HM) for functions with negative and 
monotonic Schwar\-zian derivative (decreasing in the direction of the converging sequence). For positive
Schwarzian derivative, the SNM also converges under similar monotonicity conditions, while the HM may not converge. 
We obtain conditions for the convergence of the SNM in an interval and show how this method can be 
applied for a reliable and fast solution of some problems, like the inversion of cumulative
distribution functions (gamma and beta distributions) and the inversion of elliptic integrals.

\section{The Schwarzian-Newton method: definition, geometric interpretation and convergence properties}

Let $f$ be sufficiently differentiable. Our aim is to solve $f(x)=0$. We write
\begin{equation}
\label{uno}
f''(x)+B(x)f'(x)=0, B(x)=-f''(x)/f'(x),
\end{equation}
which is a second order ODE that we can transform to normal form by setting
\begin{equation}
\label{dos}
\Phi=\exp\left(\frac12\int B\right) f=\Frac{f}{\sqrt{|f'|}}
\end{equation}
This leads to 
\begin{equation}
\label{tres}
\Phi''+\Omega \Phi=0
\end{equation}
with 
\begin{equation}
\label{cuatro}
\Omega=-\frac14 B^2 -\frac12 B'=\frac{1}{2}\left(\Frac{f'''}{f'}-\Frac{3}{2}\left(\Frac{f''}{f'}\right)^2\right)=\Frac{1}{2}\{f,x\}
\end{equation}
where $\{f,x\}$ is the Schwarzian derivative of $f$ with respect to $x$. 

As is well know, the application of Newton's method to the function $\Phi$ leads to the HM 
\cite{Brown:1977:OHV}:
\begin{equation}
\label{Halley}
x_{n+1}=g(x_n)=x_{n}-\Frac{\Phi (x_n)}{\Phi' (x_n)}=x_n-\Frac{f(x_n)}{f'(x_n)-\Frac{f''(x_n)}{2 f'(x_n)}f(x_n)}.
\end{equation}
This is a third order method. This is easy to check by considering that if $\Phi (\alpha)=0$ then 
$\Phi''(\alpha)=0$ (provided $\Omega$ is bounded at $\alpha$), and then $g'(\alpha)=g''(\alpha)=0$.

An improvement in the order of convergence 
can be obtained by integrating approximately 
the differential equation 
(\ref{tres}) (or the associated Riccati equation), similarly as in \cite{Segura:2010:RCO}. This gives a fourth order 
fixed point method. 

We start from the Riccati equation satisfied $h'(x)=1+\Omega(x) h(x)^2$, $h=\Phi/\Phi'$. Now let
$\alpha$ such that $h(\alpha)=0$ (and then $f(\alpha)=0$) and assume for the moment that $\Omega (x)>0$, we have:
\begin{equation}
\label{integ}
x-\alpha=\int_\alpha^x \Frac{h'(t)}{1+\Omega (t)h^2(t)}dt\approx \Frac{1}{\sqrt{\Omega(x)}}\arctan(\sqrt{\Omega(x)} h(x)). 
\end{equation}
where the approximation consists in taking $\Omega(x)$ constant in the integration. From this, we obtain an approximation
for $\alpha$. We can iterate these approximations and obtain the fixed point method $x_{n+1}=g(x_n)$ with

\begin{equation}
\label{atan}
g(x)=x-\Frac{1}{\sqrt{\Omega}}\arctan\left(\sqrt{\Omega} \Frac{\Phi }{\Phi'}\right),
\end{equation}
where
$$
\Frac{\Phi}{\Phi'}=\Frac{f}{\Frac{B}{2}f+f'}=\Frac{f}{f'-\Frac{f''}{2f'}f}.
$$
In other words, the method consists in computing the solutions of 
$\widetilde{h}'(x)=1+\Omega(x_n)\widetilde{h}(x)^2$ passing through the point $(x_n,h(x_n))$ and then computing $x_{n+1}$ from 
$\widetilde{h}(x_{n+1})=0$.

For the case $\Omega <0$, and using
$\arctan(\sqrt{A})/\sqrt{A}={\rm arctanh}(\sqrt{-A})/\sqrt{-A}$ when $A<0$, we can write:
\begin{equation}
\label{atanh}
g(x)=x-\Frac{1}{\sqrt{|\Omega}|}{\rm arctanh}\left(\sqrt{|\Omega|} \Frac{\Phi}{\Phi'}\right).
\end{equation}

The iteration functions (\ref{atan}) and (\ref{atanh}) define the SNM, $x_{n+1}=g(x_n)$; we observe that the SNM has the HM as limiting case for $\{f,x\}\rightarrow 0$.

\begin{remark}
The SNM is exact for functions with constant Schwarzian derivative because the approximate
integration in (\ref{integ}) becomes exact in this case. In this sense, it is the analogous to the standard Newton method,
which is exact for functions with constant ordinary derivative; this is why we call it Schwarzian-Newton method.
\end{remark}

A straightforward computation shows that the SNM satisfies: $g'(\alpha)=g''(\alpha)=g'''(\alpha)=0$ 
while $g^{(4)} (\alpha)=2\Omega' (\alpha)$ (where $f(\alpha)=0$); this implies, 
denoting the errors by $\epsilon_k =x_k-\alpha$:
$$
\epsilon_{n+1}=\Frac{\Omega '(\alpha)}{12}\epsilon_n^4 +{\mathcal O} (\epsilon_n^5).
$$
Therefore, this is a fourth order method involving third order derivatives. The same would be true, for example, for a
fourth order Newton-like method (for example of the type considered in \cite{Gil:2012:EAA}); 
the essential difference will be in the good non-local convergence properties of the SNM (see Theorem \ref{co1} 
and Corollary \ref{coco2}).

The SNM is related to the method described in \cite{Segura:2010:RCO}, which is a method for computing zeros of functions which are solution of second order differential equations. Differently, the SNM can be applied to the inversion of any function provided the Schwarzian derivative can be computed; therefore, it requires that the first three derivatives of the function are available.

As we mentioned, the SNM can be interpreted as an improvement of the HM.  
Next we study in more detail the connection of the SNM with the HM, particularly with respect to their geometrical interpretation.

\subsection{Geometrical interpretation of the Halley and Schwarzian-Newton methods}

For brevity, let us denote
\begin{equation}
\tan (\lambda,x)=\Frac{1}{\sqrt{\lambda}}\tan (\sqrt{\lambda}x)=\left\{ 
\begin{array}{crl}
\Frac{1}{\sqrt{\lambda}}\tan (\sqrt{\lambda}x)&,&\lambda>0,\\
x&,&\lambda=0,\\
\Frac{1}{\sqrt{-\lambda}}{\tanh} (\sqrt{-\lambda}x)&,&\lambda<0,
\end{array}
\right.
\end{equation}
and similarly for the inverse function $\arctan (\lambda,x)=\arctan(\sqrt{\lambda}x)/\sqrt{\lambda}$. In terms of this, the SNM reads
\begin{equation}
\label{newno}
g(x)=x-\arctan\left(\Frac{1}{2}\{f,x\},\Frac{f}{f'-\Frac{f''}{2f'}f}\right)
\end{equation}

It is easy to check that the more general functions for which the Schwarzian derivative is constant are 
\begin{equation}
\label{exasc}
h(x)=\Frac{\tan (\lambda ,x) +A}{B \tan (\lambda, x) +C},
\end{equation}
with $\{h,x\}=2 \lambda$. This can be checked by direct integration of the differential equation $\{h,x\}=2 \lambda$, 
$\lambda$ constant.  

In particular, the most general function with zero Schwarzian derivative derivative is 
\begin{equation}
\label{Haex}
h(x)=\Frac{x+A}{Bx+C},
\end{equation} 
which is the set of functions for which the HM is exact. This is consistent
with the fact that the SNM has  the HM as limiting case when $\{f,x\}\rightarrow 0$.

Precisely because the HM is exact for functions of the type (\ref{Haex}), a way to obtain the 
HM for computing the roots of $f(x)=0$ is by considering an osculating curve. We have
the following well-known result:

\begin{theorem} 
\label{HG}
Let $h(x)$ as in (\ref{Haex}) and 
define $y(x)=h(x-x_n)$, the HM (\ref{Halley}) is obtained by setting 
 $y(x_n)=f(x_n)$, $y'(x_n)=f'(x_n)$, $y''(x_n)=f''(x_n)$ and $y'''(x_n)=f'''(x_n)$ (thus determining the three
 constants) and obtaining
$x_{n+1}$ from $y(x_{n+1})=0$. The three constants are given by
\begin{equation}
A=\Frac{2 f(x_n)f'(x_n)}{D(x_n)},\,B=-\Frac{f''(x_n)}{D(x_n)},\,C=\Frac{2f'(x_n)}{D(x_n)}
\end{equation}
where
\begin{equation}
D(x_n)=2f'(x_n)^2-f(x_n)f''(x_n)
\end{equation}
\end{theorem}
Because of this result, the HM is also known as the method of tangent hyperbolas \cite{Salehov:1952:OTC}.

Similarly, because the most general function for which the SNM is exact is (\ref{exasc}), 
an alternative way of construction of the method in terms of osculating curves becomes available, as the next theorem shows.

\begin{theorem}
\label{tanhosc}
Let $h(x)$ as in (\ref{exasc}) and 
define $y(x)=h(x-x_n)$, the SNM (\ref{newno}) is obtained by setting 
 $y(x_n)=f(x_n)$, $y'(x_n)=f'(x_n)$, $y''(x_n)=f''(x_n)$ and $y'''(x_n)=f'''(x_n)$ (thus determining the four constants) and obtaining
$x_{n+1}$ from $y(x_{n+1})=0$. The constant $\lambda$ is given by 
\begin{equation}
\lambda=\Omega (x_n),\,\Omega(x)=\Frac{1}{2}\{f,x\} 
\end{equation}
and the other three constants $A$, $B$ and $C$ are as in Theorem \ref{HG}.
\end{theorem}

\noindent
{\it Proof.} The four conditions on $y(x)$ give
\begin{equation}
\label{relas}
\begin{array}{ll}
f(x_n)=\Frac{A}{C},& f'(x_n)=\Frac{C-AB}{C^2}\\
f''(x_n)=-2\Frac{B(C-AB)}{C^3}, & f'''(x_n)=\Frac{2}{C^4} (3B^2+\lambda C^2) (C-AB)
\end{array}
\end{equation}

The method is obtained by setting $y(x_{n+1})=0$, which gives
$$
x_{n+1}=g(x_n)=x_n-\arctan (\lambda, A).
$$
and this function $g(x)$ will be the same as (\ref{newno}) if
$$
\Frac{1}{A}=\Frac{f'(x_n)}{f(x_n)}-\Frac{f''(x_n)}{2f'(x_n)},
$$
and
$$
\lambda =\Frac{1}{2}\left\{\Frac{f'''(x_n)}{f'(x_n)}-\Frac{3}{2}\left(\Frac{f''(x_n)}{f'(x_n)}\right)^2\right\},
$$
which is immediate to check using (\ref{relas}).

And using these last two relations together with (\ref{relas}) we readily obtain the coefficients.
$\square$

\subsection{Non-local convergence properties of the Schwarzian-Newton method}

Before analyzing the non-local convergence properties of the SNM, we recall a result of convergence in
an interval for the HM, which is limiting case of the SNM. 

\begin{theorem}
\label{HalleyC}
Let $f$ with $f'\neq 0$ and $f'''$ continuous in an interval $J$ and let $\alpha\in J$ such
that $f(\alpha)=0$. Then, if $\{f,x\}<0$ in $J$ the HM converges monotonically 
to $\alpha$ for any starting value $x_0\in J$. 
\end{theorem}

\noindent
{\it Proof.}
Let us consider that $f'>0$ (if $f'>0$ we can proceed with the substitution $f\rightarrow -f$). As discussed before,
the HM is equivalent to the application of Newton's method to the function $\Phi=f/\sqrt{f'}$. We notice that
$$
\Phi''(x)=-\Frac{1}{2}\{f,x\}\Phi (x)
$$
and because we are assuming that $\{f,x\}<0$ then $\Phi(x)\Phi''(x)>0$ for all $x\in J\setminus\{\alpha\}$.

On the other hand, because $f'>0$ then $\Phi'>0$. Indeed, $f'>0$ implies that $\Phi(x)<0$ if $x<\alpha$ and 
$\Phi (x)>0$ if $x>\alpha$ (and the same for $\Phi''(x)>0$ because $\Phi(x)\Phi''(x)>0$ for all $x\in J\setminus\{\alpha\}$). 
Then $\Phi'(\alpha)>0$ and because $\Phi''(x)>0$ for $x>\alpha$ and $\Phi''(x)<0$ for $x<\alpha$ we have
$\Phi'(x)>0$ for all $x\in J$.

Therefore, $\Phi$ is strictly increasing in $J$ and $\Phi(x)\Phi''(x)>0$ for all $x\in J\setminus\{\alpha\}$, which
guarantees monotonic convergence of the Newton method to the zero $\alpha$ of the function $\Phi$. $\square$

\vspace*{0.3cm}
For an alternative formulation and proof of the Theorem \ref{HalleyC} see \cite{Melman:1997:GAC}. Similar results and geometrical 
proofs for other related third order methods can be found in  
\cite{Amat:2003:GCO}.

From the proof of Theorem \ref{HalleyC}, it is clear that in the case of positive Schwarzian derivative it is not
possible to guarantee convergence in an interval because the convexity properties of $\Phi$ are opposite to those
which guarantee convergence of the Newton method. 
Contrarily, we will see that conditions of non-local monotonic convergence can be found for
the SNM both for positive and negative Schwarzian derivative, but we will need an additional monotonicity
condition.

\begin{remark}
\label{rema}
From now on, we will consider that the hypotheses of theorem \ref{HalleyC} ($f'\neq 0$ and $f'''$ continuous in an interval 
$J$ and with $\alpha \in J$ such that $f(\alpha=0$)) hold in the
subsequent results.
\end{remark}

Before proving the convergence theorems for the SNM it is important to analyze the behavior of the function 
$h(x)=\Phi(x)/\Phi'(x)$ depending
on the sign and monotonicity properties of the Schwarzian derivative. 

\begin{lemma} 
\label{previ}
Under the hypotheses of Remark \ref{rema}, the following holds:
\begin{enumerate}

\item{} If $\{f,x\}>0$ $h(x)$ is strictly increasing when it is defined and it may have
one or two singularities (one smaller and one larger than $\alpha$).

\item{}If $\{f,x\}<0$, $h(x)$ is strictly increasing if $|h(x)|<|\Omega(x)|^{-1/2}$, 
$\Omega(x)=\Frac{1}{2}\{f,x\}$ and strictly decreasing if $|h(x)|>|\Omega(x)|^{-1/2}$ (when it is defined).
$h(x)$ as at most one zero and at most one singularity. 

\item{}Let $\{f,x\}<0$ and decreasing (increasing) in $J$. If $x_{-}$ is such that $h'(x_{-})<0$ then $h'(x)<0$ and 
$h(x)\neq 0$ for all
$x\in J$ greater (smaller) than $x_{-}$. 

\end{enumerate}
\end{lemma}

The proof of Lemma \ref{previ} follows easily from graphical arguments. We give the proof in the Appendix.

Using the these properties of the function we can now prove results for the convergence of the SNM in intervals.

\begin{theorem}
\label{co1}
Under the hypotheses of Remark \ref{rema}, the following holds:
\begin{enumerate}
\item{} Let $\{f,x\}$ be decreasing in $I=[a,\alpha]\subset J$.
If $\{f,x\}<0$ in $J$ then 
the SNM converges monotonically to $\alpha$ for any starting value
$x_0 \in[a,\alpha]$. If $\{f,x\}>0$ in part of the interval, the
same is true if, in addition, the SNM iteration satisfies $g(a)>a$.
\item{} Let $\{f,x\}$ be increasing in $I=[\alpha, b]\subset J$,
If $\{f,x\}<0$ in $J$ then 
the SNM converges monotonically to $\alpha$ for any starting value
$x_0 \in[\alpha, b]$. If $\{f,x\}>0$ in part of the interval, the
same is true if, in adddition, the SNM iteration satisfies $g(b)<b$.
\end{enumerate}

\end{theorem}

\noindent
{\it Proof.} We prove the case 1.a. The case 1.b is proved in a similar way.

First we have to prove that the SNM is defined for all $x\in [a,\alpha]$. For this,
we need $|\Omega (x)h(x)^2|<1$ in $[a,\alpha]$ (so that the $\mbox{arctanh}$ function in Eq. (\ref{atanh}) is defined); 
but this is necessarily so, because if there existed a 
$x_{-}$ such that $|h(x_{-})|>|\Omega(x_{-})|^{-1/2}$ then $h'(x_{-})<0$, and because $\{f,x\}<0$, Lemma \ref{previ}
guarantees that 
$h(x)\neq 0$ for all $x>x_{-}$, in contradiction with the fact that $h(\alpha)=0$, $x_{-}<\alpha$.

The monotonic convergence follows from the fact that the method consists in computing the solution of 
\begin{equation}
\label{unoh} 
\widetilde{h}'(x)=1+\Omega(x_n)\widetilde{h}(x)^2
\end{equation}
passing through the point $(x_n,h(x_n))$, where $h(x)=\Phi(x)/\Phi'(x)$ satisfies 
\begin{equation}
\label{dosh}
h'(x)=1+\Omega(x) h(x)^2, 
\end{equation}
and then computing the next iteration $x_{n+1}$ 
from $\tilde{h}(x_{n+1})=0$. 
Now, given $x_n\in [a,\alpha)$ and because $\Omega(x)<\Omega (x_n)$
for $x_n<x\le \alpha$, the graph of the solution of (\ref{unoh}) lies above the graph of $h(x)$, which is solution
of (\ref{dosh}); as a consequence the function $\widetilde{h}$ crosses the $x$-axis before $h(x)$ does. 
Therefore $x_{n+1}<\alpha$,
but also $x_n<x_{n+1}$ because $h(x_n)<0$ and this implies that $x_{n+1}=g(x_n)>x_n$. 
We have a monotonically increasing sequence bounded by $\alpha$ and converging to $\alpha$ which
is the only fixed point of $g(x)$ in $J$.

For $\Omega (x)>0$ the same proof is valid except for the fact that we must guarantee that $h(x)$ is defined for all $x$
in $(a,\alpha)$, for which we need that no value $x_{\infty}\in (a,\alpha)$ exist such that 
$\Phi'(x_{\infty})= 0$. But because $h(x)$ is monotonically increasing
for $\Omega (x)$ and there are no zeros of $h(x)$ in $(a,\alpha)$, if such value 
$x_{\infty}$ existed $h(x)=\Phi(x)/\Phi'(x)$ would change sign at $x_{\infty}$ an therefore $h(a)>0$; in that case 
$g(a)<a$. Therefore, if $g(a)>a$, $h(x)$ is defined for all $x$ in $(a,\alpha)$ (and also in $x=a$ because we are assuming
that $g(a)$ is defined).  $\square$

\vspace*{0.3cm}
As corollary of the previous theorem we have:

\begin{corollary}
\label{coco2}
Under the hypotheses of Remark \ref{rema}, if $\{f,x\}$ has one and only one extremum at $x_e\in J$ and it is a maximum, then
\begin{enumerate}
\item{}If $\{f,x\}$ is negative the SNM converges monotonically to $\alpha$ starting from $x_0=x_e$
\item{}If $(x_e-\alpha)(x_e-g(x_e))>0$ the SNM 
converges monotonically to $\alpha$ starting from $x_0=x_e$.
\end{enumerate}
\end{corollary}

Observe that the the second result
does not depend on the sign of $\{f,x\}$.

\vspace*{0.5cm}

We end this section with a comparison of the convergence properties of the HM and the SNM. We wrote the SNM $x_{n+1}=g(x_n)$ 
as $g(x)=x-\arctan(\Omega(x),h(x))$; the HM corresponds to $g(x)=x-\arctan(0,h(x))$. And because
$\arctan(\lambda,1)$ is decreasing as a function of $\lambda>-1$, we have, when $g(x)$ is real, that
$$
|\arctan(\Omega (x),h(x))|>|\arctan(0,h(x))| \mbox{ if } \Omega(x)<0
$$
and 
$$
|\arctan(\Omega (x),h(x))|<|\arctan(0,h(x))| \mbox{ if } \Omega(x)>0
$$
As a consequence
\begin{theorem}
The steps of the SNM ($x_{n+1}-x_n$) are of the same sign and greater (smaller) in absolute value than 
those for HM when $\{f,x\}$ is negative (positive).
\end{theorem}

For the case of negative $\{f,x\}$ this means that, when  Theorem \ref{co1} and Corollary \ref{coco2}
hold, the HM also converges monotonically. This is as expected on account of Theorem \ref{HalleyC}. And even more
interestingly:

\begin{corollary}
\label{menosite}
If $\{f,x\}<0$ in $J$ and $\{f,x\}$ is decreasing in $I=[a,\alpha]\subset J$ (or 
increasing in $I=[\alpha,b]\subset J$), $f(\alpha)=0$, 
the SNM converges monotonically to $\alpha$ (within a prescribed accuracy) in less iterations than the HM for any $x_0\in I$.
\end{corollary}

Contrarily, for the case of $\{f,x\}>0$ the steps of the SNM are smaller in absolute value than those of the HM. In this
case, there are no convergence results for the HM in an interval; the HM may not converge. For instance, considering the trivial
case of computing the root of $f(x)=\tan (x)$ in $(-\pi/2,\pi/2)$ the SNM is exact because $\{f,x\}$ is constant 
($\Omega (x)=1$ and $h(x)=\tan(x)$ and $g(x)=x-\arctan(\tan(x))=0$,) 
while the HM is given by $g(x)=x-\tan(x)$ and gives values outside
$(-\pi/2,\pi/2)$ if $x_0$ is close to $\pm \pi/2$.

\section{Applications} 

\subsection{Unimodal cumulative distribution functions}

We consider two examples of functions defined as 
$$
f(x)=\displaystyle\int_a^x \rho (x),x\in [a,b],
$$
with $\rho(x)>0$ and normalized in such a way that $f(b)=1$. In this case $f(x)$ is a cumulative distribution function
with probability density function $\rho (x)$. We call this unimodal if $\rho (x)$ has only one relative extremum
in $[a,b]$ and it is a maximum. This functions $f(x)$ have a sigmoidal shape with an inflection point at the maximum
of $\rho (x)$.

Let us recall that in the set of functions with constant negative Schwarzian derivative, we find functions
 with a sigmoidal behavior, as is the
case of $\tanh (x)$; and as we discussed earlier, the SNM can be interpreted in terms of osculating curves
with sigmoidal form for approximating the inversion problem. This seems adequate for these type of cumulative
distributions and provides a first hint regarding why the method can be particularly useful in these cases.  

We provide some examples of application of the method for cumulative distributions such as the central gamma and beta distributions.
As we will see, the indications that the method could have good global convergence properties will be confirmed. 
We start with the case of the central gamma distribution.

\subsubsection{The central gamma distribution}
 
As in \cite{Paris:2010:IGA} we denote
\begin{equation}
P(a,x)=\Frac{\gamma(a,x)}{\Gamma (a)},\,Q(a,x)=\Frac{\gamma(a,x)}{\Gamma (a)}
\end{equation}
where
$$
\gamma(a,x)=\displaystyle\int_0^x t^{a-1} e^{-t}dt,\, \Gamma(a,x)=\displaystyle\int_x^{+\infty} t^{a-1} e^{-t}dt.
$$
$P$ is the lower tail gamma distribution and $Q$ is the upper tail; we have $P+Q=1$, $P\in [0,1]$, $Q\in [0,1]$. $P$ is increasing,
$Q$ decreasing and they are both sigmoidal functions with inflection point at $x=a-1$. Close to this inflection point (particularly
for large $a$) the values of $P$ and $Q$ are similar.

The problem is either to invert $P(a,x)-p=0$ or $Q(a,x)-q=0$ ($q=1-p$); for obvious reasons, it is better numerically to invert 
the first equation when $p<1/2$ and the second in the other case. This is particularly important when $p$ (or $q$) is very small.
We take as function problem $f(x)=P(a,x)-p$ or $f(x)=q-Q(a,x)$ (which have the same derivatives). In both cases we have
$$
f''(x)+ B(x) f'(x)=0,\, B(x)=1+\Frac{1-a}{x}
$$

Considering the transformation (\ref{dos}) we arrive at
\begin{equation}
\Phi''(x)+\Omega (x)\Phi(x)=0,\,\Phi(x)=e^{x/2}x^{(1-a)/2}f(x),
\end{equation}
where
\begin{equation}
\label{omegai}
\Omega (x)=-\Frac{1}{4}\left(1+2\Frac{1-a}{x}+\Frac{a^2-1}{x^2}\right)=\Frac{1}{2}\left\{f,x\right\}
\end{equation}

We consider $a>0$, starting with the case $a\ge 1$. The case $a=1$ is trivial and exact for our method ($\Omega$ constant). 

We observe that for $a\ge 1$ $\Omega (x)<0$ for all $x>0$ and therefore we are in the case of
negative Schwarzian derivative. Indeed $\Omega(0^+)<0$, $\Omega(+\infty)<0$ and 
$$
\Omega' (x)=\Frac{1}{2x^3}(x(1-a)+a^2-1).
$$
Therefore the only relative extrema is at $x_m=a+1$ and it is a maximum, where $\Omega (x_m)=-(4(1+a))^{-1}<0$.

Then, the fixed point method is (\ref{atanh}) with $\Omega(x)$ given by (\ref{omegai}) and
$$
\Frac{\Phi(x)}{\Phi'(x)}=\Frac{f(x)}{\Frac{1}{2}\left(1+\Frac{1-a}{x}\right)f(x)+\Frac{e^{-x}x^{a-1}}{\Gamma (a)}}
$$
where $f(x)=P(a,x)-p$ or $f(x)=q-Q(a,x)$.

Convergence is monotonic starting from $x_0=x_m=a+1$ (see Corollary \ref{coco2}). 

It is instructive to compare the performance of the SNM with HM. We know, from Corollary 
\ref{menosite}, that the SNM will converge (within a prescribed accuracy)
in no more iterations than HM, furthermore, it has larger order of convergence. In addition, 
the SNM is exact for $a=1$ and therefore it is also extremely
efficient for values of $a$ close to $1$. 

\begin{figure}[tb]
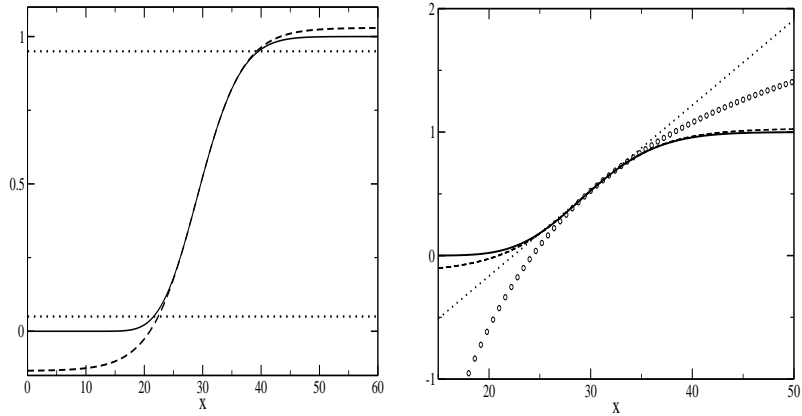

\vspace*{0.8cm}
\begin{center}
\begin{minipage}{3cm}
\centerline{\protect\hbox{\psfig{file=uno.eps,angle=0,width=5cm,height=5.4cm}}}
\end{minipage}
\hspace*{2.3cm}
\begin{minipage}{3cm}
\centerline{\protect\hbox{\psfig{file=dos.eps,angle=0,width=5cm,height=5.5cm}}}
\end{minipage}
\end{center}
\caption{Left: plot of the function $P(a,x)$ for $a=30$ (solid line), together with the osculating
curve at $x=a+1$ for the SNM method (dashed line); the horizontal lines $P=0.05$ and $P=0.95$ are 
are also shown (dotted lines). Right: plot of $P(a,x)$ for $a=30$ (solid line) and the osculating
curves at $x=a+1$ for the SNM (dashed), HM (circles) and Newton method (dotted). 
}
\label{fig1}
\end{figure}

We also provide graphical evidence of the superiority of the SNM by plotting the osculating curves at $x=x_m$ 
(see Figure \ref{fig1}); as
the figure shows, the osculating curve for the SNM is much closer to the function $P(a,x)$ than
the corresponding curves for the Newton method and the HM.

For $0<a<1$ the relative extremum $x_m=a+1$ of $\Omega (x)$ is no longer a maximum but a minimum; $\Omega (x)$
is decreasing for $x<x_m$ and increasing for $x>x_m$. We can maintain the previous iteration function but we
can no longer use as starting value $x_0=x_m$; instead, we will have monotonic convergence starting from either
large enough $x$ or small enough $x$, depending on the value of $p$. Additionally, $\Omega(x)$ becomes positive
for small enough $x$ and the resulting algorithm becomes slightly more complicated.

A way round is to consider a change of variables so that the transformed function has the same properties as before:
negative Schwarzian derivative and only one extremum (a maximum) at most. In particular, we can consider the changes
\begin{equation}
z(x)=\left\{
\begin{array}{l}
x^m/m,\, m>0\\
\log(x),\,m=0,
\end{array}
\right.
\end{equation}
as done in \cite{Deano:2004:NIF}.
In this new variable, the function $f(z)=f(x (z))$ ($x (z)$ the inverse function of $z (x)$) is such that 
$\Phi (z)=f(z)/\sqrt{\dot{f}(z)}$ satisfies $\ddot{\Phi}(z)+\Omega (z)\Phi(z)=0$ where dots mean derivative with respect
to $z$ and 
\begin{equation}
\begin{array}{l}
\Omega (z)=-\Frac{1}{4 x^{2m}}(x^2-2(a-1)x+a^2-m^2),\,\\
\dot\Omega (z)=\Frac{1}{2}\dot{x}x^{-2m-1}((m-1)x^2+(2m-1)(1-a)x+m(a^2-m^2))
\end{array}
\end{equation}
(see \cite{Deano:2004:NIF}, section 4). Two interesting cases are $m=1/2$ and $m=0$; the case $m=1$ ($z(x)=x$) was
considered before.

For $m=1/2$, $\Omega (z)$ is negative for all $x>0$ (and then $z>0$) if $a\ge 1/2$, with a maximum 
at $x=\sqrt{a^2-1/4}$ ($z=2(a^2-1/4)^{1/4)}$). The situation is then similar as before and we can perform the inversion with starting value corresponding to this maximum if $a\ge 1/2$. 

For $m=0$ ($z(x)=\log x$), $\Omega (z)$ is negative for all $x>0$ (and $z\in (-\infty,+\infty)$) if $a>0$. The only
extremum corresponds to $x=a-1$, that is $z_e=\log(a-1)$. For $a<1$, $\Omega(z)$ is strictly monotonic in ${\mathbb R}$
(because $z_e\notin {\mathbb R}$) and it is decreasing. Therefore, considering a value of $z<\log \alpha$ ($f(\alpha)=0$),
convergence is guaranteed. For $a>1$, $\Omega (z)$ is negative and with a maximum at $z=\log(a-1)$, and then the convergence
properties are similar to the case $m=1$.

Therefore, combining the case $m=1$ for $a\ge 1$ with $m=0$ for $a\in (0,1)$ ($m=1/2$ for $a\ge 1/2$ is also possible)
we have a reliable and fast method of inversion of the gamma distribution.
  
A recent algorithm for the inversion of the gamma distribution can be found in \cite{Gil:2015:GAP}. The approach was to 
compute sufficiently accurate initial approximations, depending on the range of parameters in order to ensure convergence
for a high order Newton method. The SNM provides a simpler and efficient method of computation with guaranteed convergence,
particularly for not to small $p$ or $q=1-p$; for instance, when $0.1<p<0.9$, three iterations are enough for $20$ digits
accuracy starting with $x_0=a-1$. The use of the starting values considered in \cite{Gil:2015:GAP} can, of course, 
improve the performance, particularly for small $p$ or $q$.

The noncentral gamma distribution $Q_{\mu}(x,y)$ (which becomes the central distribution $Q(\mu,y)$ for $x=0$) can also be inverted using the SNM. In particular, in \cite{Gil:2015:GAP} we applied this method for inverting $Q_{1/2}(x,y)=q$, both 
with respect to $x$ with $y$ fixed and with respect to $y$ with $x$ fixed, where
$$
Q_{1/2}(x,y)=\Frac{1}{2}\left(\mbox{erfc}(\sqrt{y}+\sqrt{x})+\mbox{erfc}(\sqrt{y}-\sqrt{x})\right).
$$

\subsubsection{The central beta distribution}

As another example of the inversion of cumulative distributions, we consider the cumulative beta distribution
\begin{equation}
I_x(a,b)=\Frac{1}{B(a,b)}\displaystyle\int_0^x t^ {a-1} (1-t)^{b-1}dt,\, x\in[0,1],\,B(a,b)=
\Frac{\Gamma (a)\Gamma (b)}{\Gamma (a+b)},
\end{equation}
with complementary function $J_x (a,b)=1-I_x (a,b)$. The problem is to invert $f(x)=I_x (a,b)-p$ (or $f(x)=q-J_x (a,b)$, $q=1-p$). 

In this case we have
\begin{equation}
\begin{array}{l}
B (x)=-\Frac{a-1}{x}+\Frac{b-1}{1-x},\\
\\
\Omega (x)=\Frac{(a-1)(b-1)}{2x(1-x)}-\Frac{1}{4}\Frac{a^2-1}{x^2}-
\Frac{1}{4}\Frac{b^2-1}{(1-x)^2}.
\end{array}
\end{equation} 

For the moment let us consider $a>1$ and $b>1$.
Because $\Omega (0^+)=-\infty$, $\Omega (1^{-})=-\infty$ and $\Omega(x)$ is differentiable in $(0,1)$, it has at least
one extremum in $(0,1)$. Now we check that when $a$ and $b$ are greater that one there is exactly one extremum and it is
is a maximum; in addition, $\Omega (x)$ is negative in $(0,1)$. This means that, similarly as in the case of the gamma
distribution, convergence of the SNM can be guaranteed by choosing as starting value the value of $x$ corresponding
to this maximum. 

That $\Omega(x)$ is negative follows from the observation that the equation $\Omega (x)=0$ does not have
real roots if $a>1$ and $b>1$. This is easy to check by writing
$$
\begin{array}{l}
\Omega (x)=-\Frac{1}{4}\Frac{P(x)}{x^2 (1-x)^2},\, P(x)=D x^2 + E x +F,\\
D =(\alpha +\beta)(\alpha +\beta +2),\,E=-(2\alpha^2+2\alpha\beta+4\alpha),\,F=\alpha(\alpha+2)
\end{array}
$$ 
where $\alpha=a-1$ and $\beta=b-1$. The discriminant is $\Delta=E^2-4DF=-2\alpha\beta (\alpha+\beta +2)<0$ and therefore there
are no real roots.

For proving there is only one relative extremum, we compute $\Omega '(x)$
$$
\begin{array}{l}
\Omega'(x)=\frac12 x^{-3} (1-x)^{-3} Q(x), \, Q(x)=G x^3 + H x^2 + I x +J,\\
G =(\alpha +\beta)(\alpha +\beta +2),
H=-3 (\alpha^2 +\alpha \beta + 2 \alpha ),\\
I=3\alpha^2 +\alpha\beta+ 6\alpha ,\,
J=-\alpha(\alpha+2)
\end{array}
$$
Then, by Descartes rule of signs we see that $Q(x)$ has either $1$ or $3$ real roots. But if it had $3$ roots then $Q'(x)$
should have two real roots. But it is easy to check that the equation $Q'(x)=0$ does not have real roots when $\alpha >0$, $\beta>0$.

The real root $x_m$ of $Q(x)$ (which gives the abscissa of the maxima of $\Omega(x)$) can be computed using standard formulas
for solving the cubic equation  $Q(x_m)=0$.

Same as for the gamma distribution, the changes of variable in \cite{Deano:2004:NIF} can be used to deal with other parameter
cases. For instance, with the change $z(x)=\log(x/(1-x))$ we arrive to
\begin{equation}
\Omega(z)=\Frac{1}{4}(-(a+b)(a+b-2)x(z)^2+2(a+b)(a-1)x(z)-a^2)
\end{equation}
which is always negative for $x\in [0,1]$ ($a$ and $b$ are positive). For $a<1$, $b<1$ $\Omega(z(x))$ has a minimum at $x_m=(a-1)/(a+b-2)$,
for $a>1$ and $b>1$ there is a maximum at $x_m$, while for the rest of cases the function is $\Omega$ is monotonic (for $x\in[0,1]$) . 
The analysis of
the monotonicity properties is more simple than in the previous case without change of variables. 
Furthermore, it appears that in some cases this alternative method is more effective \cite{prepa}.

\subsection{The incomplete elliptical integral of the second kind}

In the previous two examples we considered functions with negative Schwarzian derivative, which is the case of simpler
application of the SNM (see Theorem \ref{co1} and Corollary \ref{coco2}). But if it becomes positive, the SNM can
be also applied if the monotonicity properties of the Schwarzian derivative are available. 

As an example of this, we consider the inversion of 
\begin{equation}
f(x)=E(\sin(x),m)-pE(1,m)=\displaystyle\int_0^{x}\sqrt{1-m^2\sin^2 t}\,dt -p E(1,m)
\end{equation}
with respect to $x$, where $x\in [0,\pi/2]$, $0\le m\le 1$, $p\in [0,1]$. The function is increasing in the interval
$[0,\pi/2]$ and $f(0))=-p E(1,m)<0$, $f(\pi/2)=(1-p)E(1,m)>0$, therefore it has one and only one zero in this interval. 

The inversion of this function was
recently considered in \cite{Boyd:2012:NPC,Fukushima:2013:NIO}; we later discuss the advantages of our approach.

Proceeding as before we have
\begin{equation}
\Omega (x)=\Frac{m^2}{4}\Frac{m^2 \cos^4 x +(m^2-4)\cos^2 x +2 (1-m^2)}{(1-m^2\sin^2 x)^2}
\end{equation}

Differently to the previous cases, $\Omega (x)$ changes sign and $\Omega(x)<0$ if $x<x_c (m)$ while $\Omega(x)>0$ if $x>x_c (m)$
where
$$
x_c (m)=\arccos\left(\sqrt{\Frac{4-m^2-\sqrt{9m^4+16(1-m^2)}}{2m^2}} \right).
$$
This means we should use (\ref{atan}) when $x>x_c (m)$ and (\ref{atanh}) when $x<x_c (m)$ (or the general iteration function 
(\ref{newno})). As a function of $m$, $x_c (m)$ is increasing with $x_c (0)=\pi /4$ and $x_c (1)=\pi /2$. 

Differentiating,
\begin{equation}
\Omega '(x)=-\Frac{m^2 \sin x\cos x}{2}\Frac{3m^2(2-m^2)\cos^2 x +3m^4+m^2-4}{(1-m^2\sin^2 x)^3}
\end{equation}
Therefore we have $\Omega' (x)=0$ at $x=0$, $x=\pi/2$ and at $x_e$ if $m>2/\sqrt{7}$ with
$$
\cos^2 x_e=\Frac{3m^4 +m^2-4}{3m^4-6m^2}=1-\Frac{7m^2-4}{3m^2(2-m^2)}.
$$
It is easy to see that $x_e<x_c$ and therefore
$\Omega (x_e)<0$ and $\Omega(x)$ reaches its minimum at $x=x_e$. In the other case, when $m\le 2/\sqrt{7}$, $\Omega(x)$ is
monotonically increasing in $(0,\pi/2)$.

Therefore, considering Theorem \ref{co1} for $m\le 2/\sqrt{7}$ the SNM converges monotonically to the root, starting with 
$x_0=g(\pi/2)$, where
$$
g(\pi/2)=\Frac{\pi}{2}-\Frac{\sqrt{2(1-m^2)}}{m}\arctan\left(\Frac{m E(1,m)(1-p)}{\sqrt{2}(1-m^2)}\right).
$$
It turns out that $0<g(\pi/2)<\pi/2$ for $m\in (0,1)$, $p\in (0,1)$.

If $m>2/\sqrt{7}$, we also have monotonic convergence to the root $\alpha$ starting from $x_0=g(\pi/2)$ if $\alpha\ge x_e$; 
contrarily, if $\alpha <x_e$ we have monotonic convergence starting from $x_0=g(0)$, because $\Omega(x)$ would be
decreasing between $0$ and $\alpha$. We have 
$$
g(0)=\Frac{\sqrt{2}}{{m}}{\rm arctanh}\left(\Frac{mpE(1,m)}{\sqrt{2}}\right).
$$

For $m>2/\sqrt{\pi}$ the SNM will converge monotonically with one of the two starting values $x_0=g(0)$ or $x_0=g(\pi/2)$.
It does not seem easy to determine a priori which of the two is the correct selection, except in some cases. For
example if $g(0)>g (\pi/2)$, $g(\pi/2)$ is the best option. Even when this fact is not determined a priori, one can
build a very efficient algorithm using these starting values.

The value $g(0)$ turns out to be a good approximation for not to large values of $p$ and it is generally better than
$g(\pi/2)$. The reason for this is that $\Omega(x)$ has slower variation near $x=0$ than near
$x=\pi/2$. We have observed that considering $g(0)$ as starting value instead of $g (\pi/2)$ is the best option
when $g(0)<g(\pi/2)$ and $p$ not too large (say, $p<0.8$). 
This can be complemented with a simple approximation when $m$ is very close to $1$ (say $m>0.95$) 
using that $E(\sin x,1)=\sin x$).

Using these approximations as starting values, we have checked that a relative accuracy close to
$10^{-25}$ can be obtained with only two iterations, 
and that the accuracy is better than $10^{-40}$ in two iterations if $m$ is smaller than $0.8$. This is
clearly better than the $10^{-10}$ accuracy in three iterations of 
\cite{Boyd:2012:NPC}. Notice, in addition, that our algorithm quadruplicates
the number of exact digits in each iteration while the algorithm in \cite{Boyd:2012:NPC} 
uses Newton-Raphson, which only duplicates it. In \cite{Fukushima:2013:NIO} alternative methods are discussed which use 
accelerated bisection improved with the 
Halley method; our algorithm does not need acceleration because it is a fast high order method from the start
(faster than Halley method).

Similar ideas can be used for the inversion of other incomplete elliptic integrals. In particular, the case of the
incomplete elliptic integral of the first kind is very similar.

\section*{Appendix: proof of Lemma \ref{previ}}

\noindent
{\it Proof.}
We know that $h'(x)=1+\Omega(x)h(x)^2$ and the first two statements concerning the monotonicity are obvious. 

For the case $\Omega (x)>0$ the function is always increasing when it is defined, and necessarily the zeros
and singularities interlace. But because there is only one zero of $h(x)$ in $J$ there can be two singularities
of $h(x)$ at most.

For the case $\Omega (x)<0$, $h(x)$ is increasing
in a region symmetric around the $x$-axis ($|h(x)|<|\Omega (x)|^{-1/2}$) and decreasing outside that region. 
Then, if it has a zero $\alpha$ it is increasing at the zero, and it is easy to check by graphical arguments 
that $h(x)$ can not cross the $x$-axis again for $x>\alpha$ or $x<\alpha$. 
Similarly, proceeding with $\bar{h}=-\Phi'(x)/\Phi (x)$ and because $\Phi''+\Omega(x)\Phi (x)=0$, we have $\bar{h}'=\Omega +\bar{h}$,
which also increasing in a band around the $x$-axis. Following the same argument as before, $\bar{h}$ as one zero at most. Therefore
$h(x)$ has one singularity at most.

As for the third item, we consider the case of decreasing $\{f,x\}$; the second case is analogous. We prove that
$h (x)$ is decreasing when it is defined for $x>x_{-}$ and therefore that $h(x)\neq 0$ ($h(x)$ is necessarily increasing at its zeros, as the second statement of this theorem confirms). Because $h'(x_{-})<0$, there are two possibilities: either 
$h(x_{-})>|\Omega(x_{-})|^{-1/2}$ or $h(x_{-})<-|\Omega(x_{-})|^{-1/2}$. We start with the first case:

1. If $h(x_{-})>|\Omega(x_{-})|^{-1/2}=\lambda (x_{-})$, 
then $h(x)$ is decreasing at $x=x_{-}$ and will remain decreasing for $x>x_{-}$ with 
$h(x)>\lambda(x)$. The reason for this is that $\lambda(x)$ is decreasing in $J$ and then there can not exist
a value $x_c>x_{-}$ such that $h(x_c)=\lambda(x_c)$ (and then $h'(x_c)=0$). Indeed, because $h(x_{-})>\lambda (x_{-})$, then 
$h(x_{c}^{-})>\lambda (x_{c}^{-})$ and therefore $h'(x_{c})\le \lambda'(x_c)<0$, contradicting the fact that  $h'(x_c)=0$.

2. If $h(x_{-})<-\lambda (x_{-})$ then, because $-\lambda(x)$ is increasing, $h (x)$ will remain negative and decreasing 
($|h(x)|$ increasing) as long as it is continuous. It may 
happen that $h (x)$ has a vertical asymptote at certain $x_{\infty}>x_{-}$ such
that  $h(x_{\infty}^{-})=-\infty$ ($x_{\infty}$ would be a zero of $\Phi'(x)$). In that case, if $h(x)$ is defined for 
$x>x_{\infty}$ we would have $h (x_{\infty}^{+})=+\infty$ and $h (x)$ would remain
positive and decreasing in the rest of the interval (we are in the case discussed before). $\square$

\bibliographystyle{amsplain}
\bibliography{schwa_arxiv2}

\end{document}